\documentclass[12pt,reqno]{amsart}
\usepackage{amssymb,latexsym,amsmath,mathrsfs}

\usepackage[matrix, arrow]{xy}
\xyoption{all}

\newcommand{\A}{\mathcal{A}}

\renewcommand{\H}{\mathcal{H}}
\newcommand{\K}{\mathcal{K}}
\renewcommand{\L}{\mathcal{L}}

\newcommand{\M}{{\mathcal{M}}}
\newcommand{\C}{\mathbb{C}}
\newcommand{\N}{\mathbb{N}}
\newcommand{\T}{\mathbb{T}}
\newcommand{\Z}{\mathbb{Z}}
\newcommand{\R}{\mathbb{R}}
\newcommand{\Ff}{\mathscr{F}}
\newcommand{\Tt}{\mathfrak{T}}

\newcommand{\Od}{{\mathcal{O}_2}}

\newcommand{\ot}{\otimes}
\newcommand{\pf}{\noindent {\mbox{\textit{Proof}. }} }
\newcommand{\ie}{\textit{i.e.\/}\ } 

\newcommand{\cst}{C$^*$} 

\newtheorem{thm}{Theorem}[section]
\newtheorem{cor}[thm]{Corollary}
\newtheorem{lem}[thm]{Lemma}
\newtheorem{prop}[thm]{Proposition}

\theoremstyle{definition}
\newtheorem{defi}[thm]{Definition}
\newtheorem{exs}[thm]{Examples} 
 
\newtheorem{rem}[thm]{Remark}

\theoremstyle{remark}

\numberwithin{equation}{section}
\title{Amalgamated free products of C$^*$-bundles}

\author{Etienne Blanchard}

\date{21/11/2007}

\begin{document}



\begin{abstract} 
Given two unital continuous \cst-bundles $A$ and $B$ 
over the same compact Hausdorff base space $X$, 
we study the continuity properties of their different amalgamated 
free products over $C(X)$. 
\end{abstract} 

\maketitle
\vspace{-15pt} ${}$\hfill In memory of Gert Pedersen
\section{Introduction}
Tensor products of \cst-bundles have been much studied over the last 
decade (see for example \cite{Ri}, \cite{ElNaNe}, \cite{bla2}, \cite{KiWa}, \cite{GiMi}, \cite{Ar}, \cite{Ma}, \cite{BdWa}). 
One of the main results was obtained by Kirchberg and Wassermann 
who gave in \cite{KiWa} a characterization of the exactness 
(respectively the nuclearity) of the \cst-algebra of sections $A$ of 
a continuous bundle of \cst-algebras over a compact Hausdorff space $X$ 
through the equivalence between 
the following conditions $\alpha_e$) and $\beta_e$) 
(respectively $\alpha_n$) and $\beta_n$)): 

\smallskip 
\noindent $\alpha_e$) \textit{
The \cst-bundle $A$ is an exact \cst-algebra. 
} 

\noindent $\beta_e$) \textit{
For all continuous \cst-bundle $B$ over a compact Hausdorff space $Y$, 
the minimal \cst-tensor product $A\mathop{\ot}\limits^m B$ is a continuous \cst-bundle 
over $X\times Y$ with fibres $A_x\mathop{\ot}\limits^m B_y$. 
} 

\noindent $\alpha_n$) \textit{
The \cst-bundle $A$ is a nuclear \cst-algebra. 
} 

\noindent $\beta_n$) \textit{
For all continuous \cst-bundle $B$ over a compact Hausdorff space $Y$, 
the maximal \cst-tensor product $A\mathop{\ot}\limits^M B$ is a continuous \cst-bundle 
over $X\times Y$ with fibres $A_x\mathop{\ot}\limits^M B_y$. 
} 
\begin{rem} 
In \cite{KiWa}, the authors add to condition $\beta_e)$ 
the assumption that all the fibres $A_x$ ($x\in X$) are exact. 
But this is automatically satisfied (see \cite[Prop~3.3]{BdWa}).
\end{rem}

The case when we restrict our attention to fibrewise tensor products 
was then extensively studied in \cite{BdWa}. 
The two first assertions (respectively the two last ones) are indeed equivalent to 
the following assertion $\gamma_e$) (respectively $\gamma_n$)) introduced in \cite{bla2}, 
in case the compact Hausdorff space $X$ is perfect and second countable. 

\smallskip 
\noindent $\gamma_e$) \textit{
For all continuous $C(X)$-algebra $B$, 
the smallest completion $A\mathop{\ot}\limits^m_{C(X)} B$ 
of the algebraic tensor product $A\mathop{\odot}\limits_{C(X)} B$ amalgamated over $C(X)$ 
is a continuous \cst-bundle over $X$ with fibres $A_x \mathop{\ot}\limits^m B_x$. 
}

\smallskip 
\noindent $\gamma_n$) \textit{
For all continuous $C(X)$-algebra $B$, 
the largest completion $A\mathop{\ot}\limits^M_{C(X)}B$ 
of the algebraic tensor product $A\mathop{\odot}\limits_{C(X)} B$ amalgamated over $C(X)$ 
is a continuous \cst-bundle over $X$ with fibres $A_x \mathop{\ot}\limits^M B_x$. 
}

\smallskip 
But there are also other canonical amalgamated products over $C(X)$,
such as the completions considered 
by Pedersen (\cite{Ped}) and Voiculescu (\cite{Voi}) 
of the algebraic amalgamated free product $A\mathop{\circledast}\limits_{C(X)}B$ 
of two unital continuous \cst-bundles $A$ and $B$ 
over the same compact Hausdorff space $X$. 
The point of this paper is to  study whether analogous continuity properties hold (or not) 
for these amalgamated free products. 

\smallskip 
More precisely, we start in \S\ref{prelim} by fixing our notations and 
extending a few results available for $C(X)$-algebras 
to the framework of the operator systems which naturally appear 
when dealing with free products of $C(X)$-algebras amalgamated over $C(X)$. 
We show in \S\ref{full} that the full amalgamated free products are 
always continuous (Theorem~\ref{fullfreecont}) 
and we prove in \S\ref{red} that the exactness of  the C$^*$-algebra $A$ is sufficient 
to ensure the continuity of the reduced ones (Theorem~\ref{redfree}). 
In particular, this implies that 
any separable continuous \cst-bundle over a compact Hausdorff space $X$ 
admits a $C(X)$-linear embedding into a \cst-algebra 
with Hausdorff primitive ideal space $X$. 

\medskip 
The author would like to express his gratitude to S. Wassermann and 
N. Ozawa for helpful comments. 
He would also like to thank the referee 
for his very careful reading of several draft versions of this paper.

\section{Preliminaries}\label{prelim} 
We recall in this section a few basic definitions and constructions 
related to the theory of \cst-bundles. 

Let us first fix a few notations for operators acting on Hilbert \cst-modules (\cite[\S 13]{blac}). 
\begin{defi}\label{notationHilbertMod}
Let $B$ be a \cst-algebra and $E$ a Hilbert $B$-module. \\ 
-- For all $\zeta_1, \zeta_2\in E$, we define the rank $1$ operator $\theta_{\zeta_1, \zeta_2}$ acting on the Hilbert $B$-module $E$ by the relation \hspace{20pt} 
\begin{equation}\label{rank1} 
\theta_{\zeta_1, \zeta_2}(\xi)=\zeta_1\langle\zeta_2, \xi\rangle\quad (\xi\in E). 
\end{equation} 
-- The closed linear span of these operators is the \cst-algebra $\K_B(E)$ of \textit{compact} operators acting on the Hilbert $B$-module $E$. \\
-- The multiplier \cst-algebra of $\K_B(E)$ is (isomorphic to) the \cst-algebra $\L_B(E)$ of 
continuous adjointable $B$-linear operators acting on $E$ \\ 
-- In case $B=\C$, then $E$ is a Hilbert space and we simply denote by $\L(E)$ and $\K(E)$ 
the \cst-algebras $\L_\C(E)$ and $\K_\C(E)$. 
A basic example is the separable Hilbert space $\ell_2(\N)$ of complex valued sequences $( a_i )_{i\in\N}$ which satisfy $\| (a_i)\|^2=\sum_i |a_i|^2<\infty$. 
\end{defi}

Let $X$ be a compact Hausdorff space and 
let $C(X)$ be the \cst-algebra of continuous functions on $X$ 
with values in the complex field $\C$. 
\begin{defi} 
A $C(X)$-algebra is a \cst-algebra $A$ endowed 
with a unital $*$-homomorphism from $C(X)$ 
to the centre of the multiplier \cst-algebra $\M(A)$ of~$A$. 
\end{defi}

For all $x\in X$, we denote by $C_x(X)$ the ideal 
of functions $f\in C(X)$ satisfying $f(x)=0$. 
We denote by $A_x$ the quotient of $A$ 
by the \textit{closed} ideal $C_x(X) A$ and 
by $a_x$ the image of an element $a\in A$ in the \textit{fibre} $A_x$. 
Then the function 
\begin{equation} 
x\mapsto\| a_x\| =\inf\{\| [1-f+f(x)]a\| ,f\in C(X)\}
\end{equation}
is upper semi-continuous by construction. 
The $C(X)$-algebra is said to be \textit{continuous} 
(or to be a continuous \cst-bundle over $X$ in \cite{didou}, \cite{bla1}, \cite{KiWa}) 
if the function $x\mapsto\| a_x\|$ is actually continuous for all element $a$ in $A$. 

\begin{exs}\label{excont} 
Given a \cst-algebra $D$, the spatial tensor product $A=C(X)\ot D=C(X; D)$ 
admits a canonical structure of continuous 
$C(X)$-algebra with constant fibre $A_x\cong D$. 
Thus, if $A'$ is a \cst-subalgebra of $A$ 
stable under multiplication with $C(X)$, 
then $A'$ also defines a continuous $C(X)$-algebra. 
This is especially the case for separable exact continuous 
$C(X)$-algebras: they always admit a $C(X)$-embedding 
in the constant $C(X)$-algebra $C(X; \Od)$, 
where $\Od$ is the Cuntz \cst-algebra (\cite{bla3}). 
\end{exs} 

\medskip 
\begin{defi}\label{contfaithrep} (\cite{bla2}) 
Given a continuous $C(X)$-algebra $B$, 
a \textit{continuous field of faithful representations} 
of a $C(X)$-algebra $A$ on $B$ is a $C(X)$-linear map $\pi$ from $A$ 
to the multiplier \cst-algebra $\M(B)$ of $B$ such that, for all $x\in X$, 
the induced representation $\pi_x$ of the fibre $A_x$ in $\M(B_x)$ 
is faithful. 
\end{defi} 
Note that the existence of such a continuous field of faithful 
representations $\pi$ implies that the $C(X)$-algebra $A$ is continuous 
since the function 
\begin{equation}\label{lsc}
x\mapsto \| a_x\|=\|\pi_x(a_x)\| =\| \pi(a)_x|| =\sup\{\| (\pi (a)b)_x\| , 
b\in B \;\;\mbox{such that}\;\; \| b\|\leq 1\} 
\end{equation} 
is lower semi-continuous for all $a\in A$. 
\\ \indent 
Conversely, any \textit{separable} continuous $C(X)$-algebra $A$ admits 
a continuous field of faithful representations. 
More precisely, there always exists 
a unital positive $C(X)$-linear map $\varphi: A\to C(X)$ 
such that all the induced states $\varphi_x$ 
on the fibres $A_x$ are faithful (\cite{bla1}). 
By the Gel'fand-Naimark-Segal (GNS) construction this gives 
a continuous field of faithful representations of $A$ 
on the continuous $C(X)$-algebra of compact operators $\K_{C(X)}(E)$ 
on the Hilbert $C(X)$-module $E=L^2(A,\varphi)$. 

These constructions admit a natural extension 
to the framework of operator systems. 
Indeed, for all Banach space $V$ with a unital contractive homomorphism from $C(X)$ into the bounded linear operators on $V$, one can define the fibres $V_x=V/C_x(X) V$ and the projections $v\in V\mapsto v_x:= v+C_x(X)\, V\in V_x$ (\cite{didou}, \cite[\S 2.3]{BdKi}). 
Then the following $C(X)$-linear version of Ruan's characterization of operator spaces holds. 
\begin{prop}\label{repOpsyst} 
Let $W$ be a separable operator system 
which is a unital $C(X)$-module such that, 
for all positive integer $n$ and all $w$ in $M_n(W)$, 
the map $x\mapsto\| w_x\|$ is continuous. 
\begin{enumerate}
\item[(i)] Every unital completely positive map $\phi$ 
from a fibre $W_x$ to $M_n(\C)$ admits a $C(X)$-linear 
unital completely positive extension $\varphi:W\to M_n(C(X))$. 
\item[(ii)] There exist a Hilbert $C(X)$-module $E$ and a $C(X)$-linear 
map $\Phi: W\to\L_{C(X)}(E)$ such that for all $x\in X$, the induced map 
from $W_x$ to $\L(E_x)$ is completely isometric. 
\end{enumerate}
\end{prop}
\pf (i) Let $\zeta_n\in\C^n\ot\C^n$ be the unit vector $\zeta_n=\frac{1}{\sqrt{n}} \sum_{i=1}^n e_i\ot e_i$. 
Then the state $w\mapsto \langle\zeta_n, (id_n\ot\phi)(w)\, \zeta_n\rangle$ on $M_n(\C)\ot W_x\cong M_n(W_x)$ admits a $C(X)$-linear unital positive extension to $M_n(W)$ 
(\cite{bla1}, \cite{BdKi}). 
Thus, there is a $C(X)$-linear unital completely positive (u.c.p.) map 
$\varphi:W\to M_n(C(X))$ with 
$\varphi_x=\phi$ and $\|\varphi\|_{cb}=\|\phi\|_{cb}$ 
(\cite{Pau}, \cite{WaS}). 

\noindent (ii) 
The proof is the same as the one of Theorem 2.3.5 in \cite{EfRu2}. 
Indeed, given a point $x\in X$ and an element $w\in M_k(W)$, 
there exists, by lemma 2.3.4 and proposition 2.2.2 of \cite{EfRu2}, 
a u.c.p. map $\varphi_x$ from $W_x$ 
to $M_k(\C)$, such that 
$$\|(\imath_k\ot\varphi_x)(w_x)\|=\| w_x\|\,,$$
and one can extend $\varphi_x$ to a $C(X)$-linear u.c.p. map 
$\varphi:W\to M_k(C(X) )$ by part (i). 

For all $n\geq 1$, let $\mathfrak{s}_n$ be the set 
of completely contractive $C(X)$-linear maps from $W$ to $M_n(C(X) )$ 
and let $\mathfrak{s}=\oplus_n\,\mathfrak{s}_n$. 
Then the map $w\mapsto (\varphi(w))_{\varphi\in\mathfrak{s}}$ 
defines an appropriate $C(X)$-linear completely isometric 
representation of $W$. \qed

\begin{rem} 
Let $\{ M_n(W),\,\|\,.\,\|_n\}$ be a separable operator system such that $W$ is a continuous $C(X)$-module. 
Then the formula 

\begin{center} 
$\|w\|^\sim_n=\sup\{ \|\langle\xi\ot 1, (1_n\ot w)\, \eta\ot 1\rangle\| ; \xi, \eta\in\C^n\ot\C^n\;\mathrm{unit\; vectors} \}$
\end{center} 
for $w\in M_n(W)$ defines an operator system structure on $W$ satisfying the hypotheses of Proposition~\ref{repOpsyst}. 
\end{rem}

The Proposition \ref{repOpsyst} also induces the following $C(X)$-linear Wittstock extension: 
\begin{cor} 
Let $X$ be a compact Hausdorff space, $A$ a separable unital continuous $C(X)$-algebra 
and let $V$ be a $C(X)$-submodule of $A$. 
\\ \indent 
Then any completely contractive map $\phi$ from a fibre $V_x$ to 
$M_{k,l}(\C)$ admits a $C(X)$-linear completely contractive extension 
$\varphi: V\to M_{k,l}(C(X) )$. 
\end{cor}
\pf Let $W$ be the $C(X)$-linear operator subsystem 
$\left[\begin{array}{cc} C(X) &V\\ V^* & C(X)\end{array}\right]$ of $M_2(A)$ 
and let $\tilde\phi$ be the unital completely positive map 
from the fibre $W_x$ to $M_{k+l}(\C)$ given by 
$$\tilde\phi(\left[\begin{array}{cc}\alpha& v_x'\\
v_x^*&\beta\end{array}\right])
=\left[\begin{array}{cc}\alpha&\phi(v_x')\\ 
\phi(v_x)^*&\beta\end{array}\right]\,.$$ 
Let $\zeta=(k+l)^{-1/2}\sum_i e_i\otimes e_i\in \C^{k+l}\ot\C^{k+l}$. 
Then the associated state 
$\psi(d )=\langle\zeta, (\tilde\phi\ot\imath)(d) \zeta\rangle=
\frac{1}{k+l}\sum_{i, j}\tilde\phi_{i,j}(d_{i,j})$ 
on $M_{k+l}(W_x)$ admits a $C(X)$-linear unital 
positive extension to $M_{k+l}(W)$ (\cite{bla1}). 
Thus, there is a $C(X)$-linear completely contractive map 
$\varphi:V\to M_{k,l}(C(X))$ with 
$\varphi_x=\phi$ and $\|\varphi\|_{cb}=\|\phi\|_{cb}$ 
(\cite{Pau}, \cite{WaS}). \qed

\medskip 
We end this section with a short proof of the implication $\gamma_e)\Rightarrow\alpha_e)$ 
given in \cite{KiWa} 
\ie the characterization of the exactness of a $C(X)$-algebra $A$ by assertion~$\gamma_e$, 
if the topological space $X$ is \textit{perfect}, i.e. without any isolated point. 

\noindent $\gamma_e)\Rightarrow\alpha_e)$ 
Given two $C(X)$-algebras $A$, $B$ and a point $x\in X$, 
we have canonical $*$-epimorphisms 
$q_x:A\mathop{\ot}\limits^m B \to (A_x\mathop{\ot}\limits^m B)_x$ and 
$q_x': (A_x\mathop{\ot}\limits^m B)_x \to A_x \mathop{\ot}\limits^m B_x$. 
Further, $q_x(f\ot 1 -1\ot f)=f(x)-f(x)=0$ for all $f\in C(X)$. 
Hence $q_x$ factorizes through $A\ot_{C(X)} B$ if the $C(X)$-algebra $A$ is continuous, by \cite[proposition 3.1]{bla2}. 
If $B$ is also continuous and $A$ satisfies $\gamma_e)$, then 
$(A \mathop{\ot}\limits^m_{C(X)} B)_x \cong (A_x \mathop{\ot}\limits^m B)_x \cong A_x\mathop{\ot}\limits^m B_x$ 
and so the $C(X)$-algebra $A_x\mathop{\ot}\limits^m B$ is continuous at $x$. 
Thus, Corollary~3 of \cite{CaWa} implies that each fibre $A_x$ is exact ($x\in X$) and 
the equivalence between assertions (i) and (iv) in \cite[Thm.~4.6]{KiWa} entails 
that the \cst-algebra $A$ itself is exact.

\section{The full amalgamated free product}\label{full} 
In this section, we study the continuity of 
the full free product amalgamated over $C(X)$ 
of two unital continuous $C(X)$-algebra (\cite{Ped}, \cite{Pis}). 
By default all tensor products and free products will be over $\C$. 

\begin{defi}(\cite{VoDyNi}) 
Let $X$ be a compact Hausdorff space and 
let $A_1, A_2$ be two unital $C(X)$-algebras containing a unital copy of $C(X)$ in their centres, 
\ie $1_{A_i}\in C(X)\subset A_i$ ($i=1, 2$). \\ 
-- The \textit{algebraic free product of $A_1$ and $A_2$ with amalgamation over $C(X)$} 
is the unital quotient $A_1\mathop{\circledast}\limits_{C(X)} A_2$ 
of the algebraic free product of $A_1$ and $A_2$ over $\C$ 
by the two sided ideal generated by the differences $f1_{A_1}-f1_{A_2}$, $f\in C(X)$. \\
-- The \textit{full amalgamated free product free product} 
$A_1\mathop{\ast}\limits^f_{C(X)} A_2$ is the universal unital enveloping \cst-algebra 
of the $\ast$-algebra $A_1\mathop{\circledast}\limits_{C(X)} A_2$. \\ 
-- Any pair $(\sigma_1, \sigma_2)$ of unital $\ast$-representations of $A_1, A_2$  
that coincide on their restrictions to $C(X)$ defines 
a unital $\ast$-representation $\sigma_1\ast\sigma_2$ of $A_1\mathop{\circledast}\limits_{C(X)} A_2$, 
the restriction of which to $A_i$ coincides with $\sigma_i$ ($i=1, 2$). \\
\end{defi} 
In particular, the two unital central copies of $C(X)$ in $A_1$ and $A_2$ coherently define 
a structure of $C(X)$-algebra on $A_1\mathop{\ast}\limits^f_{C(X)} A_2$ 
and by universality, we have: 
\begin{equation}\label{fibrefull}
\forall\,x\in X\,,\quad (A_1\mathop{\ast}\limits^f_{C(X)} A_2)_x\cong 
(A_1)_x\mathop{\ast}\limits^f_{\C} (A_2)_x\,.
\end{equation} 

\begin{rem} If we fix unital positive $C(X)$-linear maps $\varphi_i:A_i\to C(X)$ and 
we set  $A_i^\circ=\ker\varphi_i$ for $i=1, 2$, 
then 
the algebraic amalgamated free product $A_1\mathop{\circledast}\limits_{C(X)} A_2$ is 
(isomorphic to) the $C(X)$-module 
$C(X)\oplus \mathop{\oplus}\limits_{n\geq 1}\mathop{\bigoplus}\limits_{i_1\not=\ldots\not= i_n}
A_{i_1}^\circ\mathop{\ot}\limits_{C(X)}\ldots\mathop{\ot}\limits_{C(X)}A_{i_n}^\circ$,
which is a $*$-algebra for the product $v.w:=v\otimes w$ and the involution $(v.w)^*=w^*.v^*$. 
\end{rem} 

Assume now that $A_1$ and $A_2$ are continuous $C(X)$-algebras. 
Then $A_1\mathop{\ast}\limits^f_{C(X)} A_2$ is also a continuous $C(X)$-algebra 
as soon as both $A_1$ and $A_2$ are separable exact \cst-algebras 
thanks to an embedding property due to Pedersen (\cite{Ped}): 
\begin{prop}\label{fullfreeexact} 
Let $X$ be a compact Hausdorff space and 
$A_1, A_2$ two separable unital continuous $C(X)$-algebras which are exact \cst-algebras. 
\\ \indent 
Then the full amalgamated free product 
$A_1\mathop{\ast}\limits^f_{C(X)} A_2$ is a continuous $C(X)$-algebra. 
\end{prop}
\pf For $i=1, 2$, let $\pi_i$ be a $C(X)$-linear embedding of $A_i$ into $C(X; \Od)$ (\S \ref{excont}). 
Then the induced $C(X)$-linear morphism $\pi_1\ast\pi_2$ from 
$A_1\mathop{\ast}\limits^f_{C(X)} A_2$ to 
the continuous $C(X)$-algebra 
$C(X; \Od)\mathop{\ast}\limits^f_{C(X)} C(X; \Od)=C(X; \Od\mathop{\ast}\limits^f\Od)$, 
is injective by \cite[Thm.~4.2]{Ped}. 
\qed

This continuity property actually always holds (Theorem~\ref{fullfreecont}). 
In order to prove it, let us first state the following Lemma 
which will enable us to reduce the problem to the separable case. 
(Its proof is the same as in \cite[2.4.7]{BdKi}.) 
\begin{lem}\label{freesepreduction} 
Let $X$ be a compact Hausdorff space, $A_1, A_2$ two unital $C(X)$-algebras and $a$ an element 
of the algebraic amalgamated free product $A_1\mathop{\circledast}\limits_{C(X)}A_2$. 
Then there exist a second countable compact Hausdorff space $Y$ s.t. $1_{C(X)}\in C(Y)\subset C(X)$ 
and separable $C(Y)$-algebras $D_1\subset A_1$ and $D_2\subset A_2$ s.t. 
\hbox{$a$ belongs to the $\ast$-subalgebra $D_1 \mathop{\circledast}\limits_{C(Y)}D_2$. }
\end{lem}

Let $e_1, e_2, \ldots$ be an orthonormal basis of $\ell^2(\N)$ and 
set $e_{i, j}:=\theta_{e_i, e_j}$ for all $i, j$ in $\N^*:=\N\setminus\{ 0\}$ (see (\ref{rank1}))
Note that $e_{i, j}$ is the rank $1$ partial isometry such that $e_{i, j} e_j=e_i$. 
Then the following critical Lemma holds: 
\begin{lem}\label{fullfreetech} 
Let $Y$ be a second countable compact space,
$D_1, D_2$ two separable unital continuous $C(Y)$-algebras 
and set $\mathcal{D}=D_1\mathop{\ast}\limits^f_{C(Y)} D_2$. 
\\ \indent 
If the element $d$ belongs to the algebraic amalgamated free product 
$D_1\mathop{\circledast}\limits_{C(Y)} D_2$, 
then the function $y\mapsto \| d_y\|_{\mathcal{D}_y}$ is continuous. 
\end{lem} 
\pf 
The map $y\mapsto \| d_y\|$ is always upper semicontinuous by (\ref{fibrefull}). 
So, it only remains to prove that it is also lower semicontinuous 
if both the $C(Y)$-algebras $D_1$ and $D_2$ are continuous. 

\medskip 
Now, any element $d$ in the algebraic amalgamated free product of $D_1$ and $D_2$ 
admits by construction (at least) one finite decomposition (that we fix) in $M_m(\C)\ot\mathcal{D}$ 
\begin{equation} e_{1,1}\ot d= d(1)\ldots d(2n)\end{equation} 
for suitable integers $m, n\in\N$ and elements $d(k)\in M_m(\C)\ot D_{\iota_k}$ ($1\leq k\leq 2n$), 
where $\iota_k=1$ if $k$ is odd and $\iota_k=2$ otherwise. 
\\ \indent 
Given a point $y\in Y$ and a constant $\varepsilon>0$, there exist 
unital $*$-representations $\Theta_1, \Theta_2$ of the fibres $(D_1)_y, (D_2)_y$ on $\ell^2(\N)$ 
and unit vectors $\xi, \xi'$ in $e_1\ot\ell^2(\N) \subset \C^m\ot\ell^2(\N)$ s.t. 
\begin{equation} 
\| d_y\|-\varepsilon < 
\Bigl|\langle\xi', e_{1,1}\ot(\Theta_1\ast\Theta_2)(d_y) \xi\rangle\Bigr| \,.
\end{equation} 
As the sequence of projections $p_k=\sum_{i=0}^k e_{i,i}$ in $\K(\ell^2(\N))$ 
satisfies $\mathop{\lim}\limits_{k\to\infty}\| (1-p_k)\zeta\|=0$ for all $\zeta\in\ell^2(\N)$, 
a finite induction implies that there is an integer $l\in\N$ such that 
$ \| (1\ot p_l)\xi\|\not=0$, $ \| (1\ot p_l)\xi'\|\not=0$ and 
the two u.c.p. maps $\phi_i(\, .\, )=p_l\Theta_i(\, .\, ) p_l$ on the fibres $(D_i)_y$ ($i=1, 2$) 
satisfy: 
\begin{equation}\label{weak} 
\bigl|\langle\xi_l', 
\Bigl[e_{1,1}\ot (\Theta_1\ast\Theta_2)(d_y)- 
(id\ot\phi_{\iota_1})(d(1)_y)\ldots(id\ot\phi_{\iota_{2n}})(d(2n)_y)\Bigr] \xi_l\rangle\bigr| < \varepsilon 
\end{equation} 
where $\xi_l=(1\ot p_l)\xi/\| (1\ot p_l)\xi\|$ and $\xi_l'=(1\ot p_l)\xi/\| (1\ot p_l)\xi'\|$ are unit vectors 
in $\C^m\ot\C^{l+1}$ which are arbitrarily close to $\xi$ and $\xi'$, respectively, 
for sufficiently large $l$. 

Let $\zeta_l\in\C^l\ot\C^l$ be the unit vector $\zeta_l=\frac{1}{\sqrt{l}} \sum_{1\leq k\leq l} e_k\ot e_k$. 
For each $i$, 
the state $e\mapsto\langle\zeta_l, (id\ot\phi_i)(e)\,\zeta_l\rangle$ on $M_{l}(\C)\ot (D_i)_y$
associated to $\phi_i$ admits a unital $C(Y)$-linear positive extension 
$\Psi_i: M_{l}(\C)\ot D_i\to C(Y)$ (\cite{bla1}). 
If $(\mathcal{H}_i, \eta_i, \sigma_i)$ is the associated GNS-Kasparov construction, 
then every $d_i\in D_i$ satisfies 
\begin{equation} 
\langle 1_l\ot\eta_i, (id\ot\sigma_i)(\theta_{\zeta_l, \zeta_l}\ot d_i) 1_l\ot\eta_i\rangle(y) 
=(id\ot\Psi_i)(\theta_{\zeta_l, \zeta_l}\ot d_i)(y) 
=\phi_i((d_i)_y) 
\end{equation} 
Let $\sigma=\sigma_1\ast \sigma_2$ be the $*$-representation of 
the full amalgamated free product $\mathcal{D}$ 
on the amalgamated pointed free product $C(Y)$-module $(\H, \eta)=\ast_{C(Y)}(\H_i, \eta_i)$ 
(\cite{VoDyNi}). 
Then
$$ {}\hspace{-3pt}
\begin{array}{rl} 
\bigl| \langle \xi_l'\ot\eta, e_{1,1}\ot \sigma(d)\, \xi_l\ot\eta\rangle\bigr|(y)\hspace{-5pt}&
=\bigl| \langle \xi_l'\ot\eta, (id\ot\sigma)(d(1))\ldots (id\ot\sigma)(d(2n))\xi_l\ot\eta\rangle\bigr|(y)\\
&=\bigl| \langle \xi_l', (id\ot\phi_{\iota_1})( d(1)_y)\ldots(id\ot\phi_{\iota_{2n}})( d(2n)_y) \xi_l \rangle\bigr|\\
&> \| d_y\|-2\varepsilon
\end{array} 
$$ 
And so, 
$\| d_y\|-2\varepsilon < 
\Bigl| \langle \xi_l'\ot\eta, e_{1,1}\ot \sigma(d)\, (1\ot p_l)\xi_l\ot\eta\rangle\Bigr|(z)
\leq \| d_z\|$ 
for all point $z$ in an open neighbourhood of $y$ in $Y$ by continuity. 
\qed

\begin{rem} 
The referee pointed out that the inequality (\ref{weak}) cannot be replaced by a norm inequality like 
$\bigl\| e_{1,1}\ot (\Theta_1\ast\Theta_2)(d_y)- (id\ot\phi_{\iota_1})(d(1)_y)\ldots(id\ot\phi_{\iota_{2l}})(d(2l)_y)\bigr\| < \varepsilon'$. \\ 
Indeed, if for instance $p\in A=M_2(\C)$ is the projection 
$p=\frac{1}{\sqrt{2}} \left[\begin{array}{cc}1&1\\1&1\end{array}\right] $, then 
$e_{1,1}\,.\,e_{2,2}=0$ \quad but\quad $p\,e_{1,1}\,p\,e_{2,2}\,p= 2^{-3/2} p\not=0\,$. 
\end{rem}

\smallskip 
\begin{thm}\label{fullfreecont} Let $X$ be a compact Hausdorff space and 
let $A_1, A_2$ be two unital continuous $C(X)$-algebras. 
Then the full amalgamated free product $\A =A_1\mathop{\ast}\limits^f_{C(X)} A_2$ 
is a continuous $C(X)$-algebra with fibres $\A_x= (A_1)_x\mathop{\ast}\limits^f (A_2)_x$ ($x\in X$). 
\end{thm}
\pf 
The $C(X)$-algebra $\A$ has fibre $\A_x=(A_1)_x \mathop{\ast}\limits^f (A_2)_x$ at $x\in~X$ 
by (\ref{fibrefull}). 
Hence it is enough to prove that for all $a$ in the dense algebraic amalgamated free product 
$A_1\mathop{\circledast}\limits_{C(X)}A_2\subset \A$, the map $x\mapsto\| a_x\|$ is lower semi-continuous. 

Let $a$ be such an element and 
choose a finite decomposition $e_{1,1}\ot a=a_1\ldots a_{2n}\in M_m(\C)\ot\A$, 
where $m, n\in\N$ and $a_k$ belongs to $M_m(\C)\ot A_1$ or $M_m(\C)\ot A_2$ according to the parity of $k$. 
By Lemma \ref{freesepreduction}, there exist 
a separable unital \cst-subalgebra $C(Y)\subset C(X)$ containing the unit $1_{C(X)}$ of $C(X)$ and 
separable unital \cst-subalgebras $D_1\subset A_1$ and $D_2\subset A_2$ such that 
each $D_i$ is a continuous $C(Y)$-algebra and 
all the $a_k$ belong to $M_m(\C)\ot D_1$ or $M_m(\C)\ot D_2$ according to the parity of $k$. 
And so $a$ also belongs to the full free product 
$\mathcal{D}=D_1\mathop{\ast}\limits^f_{C(Y)} D_2$ 
which admits a $C(Y)$-linear embedding in 
$A_1\mathop{\ast}\limits^f _{C(X)} A_2$ (\cite[Thm.~4.2]{Ped}). 
Hence it is enough to prove that the map $y\in Y\mapsto\| a_y\|_{\mathcal{D}_y}$ is also 
lower semicontinuous. 
But this follows from Lemma~\ref{fullfreetech}. 
\qed

\section{The reduced amalgamated free product}\label{red} 
Let us now study the continuity properties of 
certain reduced amalgamated free product over $C(X)$ 
of two unital continuous $C(X)$-algebras (\cite{Voi}, \cite{VoDyNi}). 
\\ \indent The main result of this section is the following: 
\begin{thm}\label{redfree} 
Let $X$ be a compact Hausdorff space and 
let $A_1, A_2$ be two unital continuous $C(X)$-algebras. 
For $i=1,2$, let $\phi_i: A_i\to C(X)$ be a unital projection 
such that for all $x\in X$, the induced state $(\phi_i)_x$ on the 
fibre $(A_i)_x$ has faithful GNS representation. 

If the \cst-algebra $A_1$ is exact, 
then the reduced amalgamated free product 
$$(A,\phi)=(A_1,\phi_1)\mathop{\ast}\limits_{C(X)} (A_2,\phi_2)$$ is 
a continuous $C(X)$-algebra 
with fibres $(A_x,\phi_x)=((A_1)_x,(\phi_1)_x) \ast \,((A_2)_x,(\phi_2)_x)$. 
\end{thm} 

The proof is similar to the one used by Dykema and Shlyakhtenko 
in \cite[\S 3]{DySh} 
to prove that a reduced free product of exact \cst-algebras is exact. 
We shall accordingly omit details except where our proof deviates from theirs. 
\begin{lem}\label{extension} 
Let $A$ be a $C(X)$-algebra and $J\triangleleft A$ be a closed two sided ideal in $A$. 
If the two $C(X)$-algebras $J$ and $A/J$ are continuous, then $A$ is also continuous. 
\end{lem} 
\pf The canonical $C(X)$-linear representation $\pi$ of $A$ on $J\oplus A/J$ is a continuous field of faithful representations. 
Indeed, if $a\in A$ satisfies $\pi_x(a_x)=0$ for some $x\in X$, 
then $(a a'+J)_x=0$ for all $a'\in A$, hence $(a+J)_x=0$, \ie $a_x\in J_x$. 
Now $a_x h_x=(a h)_x=0$ for all $h\in J$ and so $a_x=0$. \qed 

\begin{rem} 
The continuity of the $C(X)$-algebra $A$ does not imply the continuity of the quotient $A/J$. 
In fact, any $C(X)$-algebra $B$ is the quotient 
of the constant $C(X)$-algebra $A=C(X; B)=C(X)\otimes B$ 
by the (closed) two sided ideal $C_\Delta .A$, where $C_\Delta\subset C(X\times X)$ is the ideal of functions $f$ which satisfy $f(x, x)=0$ for all $x\in X$. 
\end{rem} 

\begin{lem}\label{stabeq} 
Let $B$ be a unital $C(X)$-algebra and $E$ a \textit{full} countably generated Hilbert $B$-module. 
Then $B$ is a continuous $C(X)$-algebra if and only if 
the $C(X)$-algebra $\K_B(E)$ of compact operators acting on $E$ (Definition~\ref{notationHilbertMod}) 
is continuous. 
\end{lem} 
\pf The \cst-algebra $B$ and $\K_B(E)$ are stably isomorphic by Kasparov stabilisation theorem 
(\cite[Thm. 13.6.2]{blac}), 
\ie there is a $B$-linear isomorphism \\
\centerline{$\K(\ell^2(\N) )\ot B\cong \K_B(\ell^2(\N)\ot E))\cong \K(\ell^2(\N) )\ot \K_B(E)$ 
\qquad (\cite[Ex.~13.7.1]{blac}). }
Note that this isomorphism are also $C(X)$-linear since $1_B\in C(X)\subset B$. 
As the \cst-algebra $\K(\ell^2(\N) )$ is nuclear, 
the Theorem 3.2 of \cite{KiWa} implies the equivalence between the continuity of the $C(X)$-algebras $B$, $\K(\ell^2(\N) )\ot B$ and $\K_B(E)$. \qed

\medskip 
Given a \cst-algebra $B$ and a Hilbert $B$-bimodule $E$, 
recall that the full Fock Hilbert $B$-bimodule associated to $E$ is the sum 
$\Ff_B(E)=B\oplus E\oplus (E\otimes_B E)\oplus\ldots=
\bigoplus_{n\in\N}\;E^{(\otimes_B)n}$ 
and that for all $\xi\in E$, the creation operator 
$\ell(\xi)\in\L_B(\Ff_B( E))$ is defined by 
\begin{equation}\label{Fockspace}
\begin{array}{lcll}
\bullet\quad \ell(\xi) b=\xi b&\textrm{for}&b\in B=:E^0&\mathrm{ and}\\ 
\bullet\quad  \ell(\xi)\,(\zeta_1\otimes\ldots\otimes\zeta_k)=\xi\otimes\zeta_1\otimes\ldots\otimes\zeta_k &
\textrm{for}&\zeta_1,\ldots, \zeta_k\in E\,.& 
\end{array}
\end{equation} 
Then then Toeplitz \cst-algebra $\Tt_B ( E)$ of the Hilbert $B$-module $E$ 
(called \textit{extended Cuntz-Pimsner algebra} in \cite{DySh}) is the \cst-subalgebra 
of $\L_B(\Ff_B( E))$ generated by the operators $\ell(\xi),\,\xi\in E$ (\cite{Pim}).

\begin{lem}\label{lemToeplitz} 
Let $B$ be a unital $C(X)$-algebra and let $E$ be a countably generated Hilbert $B$-bimodule 
such that the left module map $B\to\L_B (E)$
 is injective and which satisfies 
\begin{equation} \label{diamond}
f.\zeta=\zeta . f\quad\mathrm{for}\;\mathrm{all}\; \zeta\in E\;\mathrm{and}\; f\in C(X)\,.
\end{equation} 

Then the Toeplitz $C(X)$-algebra $\Tt_B (E)$ of $E$ 
is a continuous $C(X)$-algebra with fibres $\Tt_{B_x} (E_x)$ if and only if 
the $C(X)$-algebra $B$ is continuous. 
\end{lem} 

Under the assumption of this Lemma, 
the canonical $*$-monomorphism $B \to \Pi_{x\in X} B_x$ induces 
for any Hilbert $B$-module $F$ a $*$-homomorphism $a\mapsto a\ot 1$ 
from the \cst-algebra $\K_B\left(F \right)$ of compact operators acting on the Hilbert module $F$ 
to the tensor product $\K_B(F)\ot_B \left(\Pi_{x\in X} B_x\right)\cong\Pi_{x\in X} \, \K_{B_x}\left( F\ot_B B_x\right)$. 
And this map is injective as soon as there is a $B$-linear decomposition $F\cong B\oplus F'$ 
for some Hilbert $B$-module $F'$. 
After passing to the multiplier \cst-algebras, this gives for $F=\Ff_B(E)$ a $*$-monomorphism 
$\Theta=(\Theta_x)$: 
$$\L_B(\Ff_B(E)) = \M\bigl(\K_B(\Ff_B(E)) \bigr)\hookrightarrow 
\mathop\Pi\limits_{x\in X} \M\bigl(\K_{B_x}\left(\Ff_{B_x}(E_x)\right)\bigr)=
\mathop\Pi\limits_{x\in X} \L_{B_x}\left(\Ff_{B_x}(E_x)\right),$$ 
where $E_x$ is the Hilbert $B_x$-module $E_x=E\ot_B B_x\cong E/C_x(X)E$ for all $x\in X$. 

Note that for all $\xi\in E$ and $x\in X$, we have $\ell(\xi)\, C_x(X)\Ff_B(E)\subset C_x(X)\Ff_B(E)$ 
by (\ref{diamond}) and so the element $\Theta_x(\ell(\xi) )$ satisfies 
the same creation rules (\ref{Fockspace}) as the creation operator $\ell(\xi_x)$, 
where $\xi_x=\xi\ot 1_{B_x} \in E_x$. 
So, the restriction of $\Theta$ to $\Tt_B(E)$ takes values 
in the product $\mathop\Pi\limits_{x\in X} \Tt_{B_x}(E_x)$.

\medskip 
\noindent \textit{Proof of Lemma~\ref{lemToeplitz}.}
The continuity of the $C(X)$-algebra $\Tt_B ( E)$ clearly implies the continuity 
of the $C(X)$-algebra~$B$ since $B$ embeds $C(X)$-linearly in $\Tt_B ( E)$. 

Suppose conversely that the $C(X)$-algebra $B$ is continuous. 
Let $\widetilde{E}$ be the full countably generated Hilbert $B$-bimodule $\widetilde{ E}= E\oplus B$. 
Then the $C(X)$-algebra $\K_B(\Ff_B(\widetilde{E}) )$ of compact operators acting on the Hilbert $B$-module $\Ff_B(\widetilde{E})$ is a continuous $C(X)$-algebra by Lemma~\ref{stabeq}. 
Hence it is enough to prove that the Toeplitz \cst-algebra $\Tt_B (\widetilde{ E})$ 
admits a continuous field of faithful representations 
$\Tt_B (\widetilde{ E})\to\L_B(\Ff_B(\widetilde{E}) )=\M( \K_B(\Ff_B(\widetilde{E})) )$ 
since $\Tt_B ( E)$ embeds in $\Tt_B (\widetilde{ E})=\Tt_B( E\oplus B)$ 
(\cite{Pim} or \cite[\S 3]{DySh}). 

\noindent\textit{Step 1.} \textit{
Let $\beta$ be the action of the group $\T=\R/\Z$ 
on $\Tt_B (\widetilde{ E})$ 
determined by $\beta_t(\ell(\zeta))=e^{2i\pi t}\ell(\zeta)$. 
Then the fixed point \cst-subalgebra $A$ under this action is a continuous 
$C(X)$-algebra and 
the map $A_x\to \Tt_{B_x} (\widetilde{E}_x)$ is injective for all $x\in X$. }

Define the increasing sequence of $B$-subalgebra 
$A_n\subset A$ generated by the words of the form 
$w=\ell(\zeta_1)\ldots \ell(\zeta_k)\,\ell(\zeta_{k+1})^*\ldots\ell(\zeta_{2k})^*$ 
with $k\leq n$. Let also $A_0=B$. As $A=\overline{\cup A_n}$, it is 
enough to prove that each $A_n$ is continuous with appropriate fibres. 

Let $\widetilde{ E}_0=B$ and $\widetilde{ E}_n=\widetilde{ E}\ot_B\ldots\ot_B\widetilde{ E}=
\widetilde{ E}^{(\ot_B) n}$ for $n\geq 1$. 
Define also the projection $P_ n\in \L_B(\Ff_B(\widetilde{E}) )$ 
on the Hilbert $B$-bimodule $F_n=\oplus_{0\leq k\leq n} \widetilde{E}_k$. 
In the decomposition $\Ff_B(\widetilde{E} )\cong F_n\ot_B\Ff_B(\widetilde{E}^{\otimes_B(n+1)})$, 
$A_n$ acts on $\Ff_B(\widetilde{E})$ as $\Theta_n(A_n)\ot 1$, where $\Theta_n(a)=P_n aP_n$ for $a\in A_n$. 
Thus the map $\Theta_n$ is faithful on $A_n$. 
Further, the kernel of the restriction map $F_n=F_{n-1}\oplus \widetilde{E}_n\to F_{n-1}$ is the $B$-module generated by the words 
$w=\ell(\zeta_1)\ldots \ell(\zeta_n)\,\ell(\zeta_{n+1})^*\ldots\ell(\zeta_{2n})^*$ of length $2n$, 
which is isomorphic to $\K_B(\widetilde{E}_n)$. 
Hence, we have by induction $C(X)$-linear split exact sequences 
\begin{equation} 
0\to \K_B(\widetilde{ E}_n)\to A_n\to A_{n-1}\to 0\,,
\end{equation} 
and so, each $C(X)$-algebra $A_n$ is continuous by lemma \ref{extension}. 

\smallskip\noindent\textit{Step 2.} \textit{
The Toeplitz \cst-algebra $\Tt_B (\widetilde{ E})$ 
is isomorphic to the crossed product $A\rtimes_\alpha \N$, 
where $\alpha: A\to A$ is the injective $C(X)$-linear endomorphism 
$\alpha(a)=LaL^*$, with $L=\ell(0\oplus 1_B)$ 
(\cite[Claim 2.1.3]{DySh}). 
Hence $\Tt_B (\widetilde{ E})$ is a continuous field with fibres 
$\bigl(\Tt_B (\widetilde{ E})\bigr)_x\cong A_x\rtimes\N\cong \Tt_{B_x} (\widetilde{ E}_x)$ 
for $x\in X$. }

The \cst-algebra $\Tt_B(\widetilde{E})$ is generated by $A$ and $L$. 
Hence it is isomorphic to the $C(X)$-algebra $A\rtimes_\alpha\N$ (\cite[claim 3.4]{DySh}). 
Let us now study the continuity question. 

Let $\overline{A}$ be the inductive limit of the system 
$A\mathop\to\limits^\alpha A\mathop\to\limits^\alpha \ldots$ with corresponding $C(X)$-linear monomorphisms $\mu_n: A\to \overline{A}$ ($n\in\N$). 
It is a continuous $C(X)$-algebra since 
$\bigcup\mu_n(A)$ is dense in $\overline{A}$ and 
the map $x\in X\mapsto \| \mu_n(a)_x\|=\| a_x\|$ is continuous for all $(a, n)\in A\times\N$. 
Let $\overline{\alpha}:\overline{A}\to \overline{A}$ be the $C(X)$-linear automorphism given by 
$\overline{\alpha}\bigl(\mu_n(a) \bigr) = \mu_n\bigl( \alpha(a) \bigr)$, 
with inverse $\mu_n(a)\mapsto \mu_{n+1}(a)$. 
Then the crossed product $\overline{A}\rtimes_{\overline{\alpha} }\Z$ is continuous over $X$ 
since $\Z$ is amenable (\cite{Ri}). 
Hence, if $p\in \overline{A}$ is the projection $p=\mu_0(1_A)$, the hereditary $C(X)$-subalgebra 
$p\left( \overline{A}\rtimes_{\overline{\alpha} }\Z \right) p =: A\rtimes_\alpha\N$ (\cite{DySh}) is also continuous with fibres $A_x\rtimes_{\alpha_x}\N$ ($x\in X$). 
\qed

\noindent \textit{Proof of Theorem~\ref{redfree}.} 
By density, it is enough to study the case of elements $a$ in the algebraic amalgamated free product $A_1\mathop{\circledast}\limits_{C(X)} A_2$. 
But, for any such $a$, there are separable unital \cst-subalgebras $C(Y)\subset C(X)$, $D_1\subset A_1$ and $D_2\subset A_2$ with same units such that 
$a$ also belongs to $D_1\mathop{\circledast}\limits_{C(Y)} D_2$ by Lemma \ref{freesepreduction}. 
And the reduced free product $(D_1,\phi_1)\mathop{\ast}\limits_{C(Y)} (D_2,\phi_2)$ 
embeds $C(X)$-linearly in $(A,\phi)=(A_1,\phi_1)\mathop{\ast}\limits_{C(X)} (A_2,\phi_2)$ 
by \cite[Thm. 1.3]{BlDy}. 
Further, any \cst-subalgebra of a exact \cst-algebra is exact. 
Thus, one can assume in the sequel 
that the compact Hausdorff space $X$ is second countable and 
that the $C(X)$-algebras $A_1, A_2$ are separable \cst-algebras.

\bigskip
If the \cst-algebra $A_1$ is exact, then 
the $C(X)$-algebra $B=A_1\mathop{\otimes}\limits_{C(X)} A_2$ is continuous 
with fibres $B_x=(A_1)_x\mathop{\ot}\limits^m (A_2)_x$ by $\gamma_e$), and 
the conditional expectation $\rho=\phi_1\otimes\phi_2:B\to C(X)$ 
is a continuous field of states on $B$ such that each $\rho_x$ has faithful GNS representation 
($x\in X$).

Let $E$ be the full countably generated Hilbert $B$,$B$--bimodule $L^2(B,\rho)\ot_{C(X)} B$ and 
let $\Ff_B( E)=B\oplus \left( L^2(B,\rho)\mathop{\otimes}\limits_{C(X)}B \right) \oplus 
\left( L^2(B,\rho)\mathop{\otimes}\limits_{C(X)} L^2(B,\rho)\mathop{\otimes}\limits_{C(X)} B\right) \oplus
\ldots$ be its full Fock bimodule. 
Let also $\xi=\Lambda_{\phi}(1)\ot 1\in E$. 
As observed in \cite[Claim 3.3]{DySh}, the Toeplitz \cst-algebra $\Tt_B(E)\subset \L_B(\Ff_B(E) )$ is 
generated by the left action of $B$ on $\Ff_B(E)$ and the operator $\ell(\xi)$, 
because $\ell(b_1\xi b_2)=b_1\ell(\xi) b_2$ for all $b_1, b_2$ in $B$. 

Consider the conditional expectation $\mathfrak{E}: \Tt_B ( E)\to~B$ defined by compression with the orthogonal projection from $\Ff_B( E)$ into the first summand $B\subset\Ff_B(E)$. 
Then Theorem~2.3 of \cite{Sh} implies that 
$B$ and the $C(X)$-algebra generated by the non trivial isometry $\ell(\xi)$ are free with amalgamation over $C(X)$ in $(\Tt_B ( E),\rho\circ\mathfrak{E})$ because $\ell(\xi)^*\,b\,\ell(\xi)= \rho(b)$ for all $b\in B$. 

By \cite{Sh}, the restriction of $\mathfrak{E}$ to the \cst-subalgebra $C^*(\ell(\xi) )\subset \Tt_B ( E)$ 
takes values in $C(X)$ and 
there exists a unitary $u\in C^*(\ell(\xi) )$ s.t. $\mathfrak{E}(u^k)=0$ for every non-zero integer $k$. 
The two embeddings $\pi_i: A_i\to \Tt_B(E)$ ($i=1, 2$) given by 
$\pi_1(a_1)=u\, (a_1\ot 1)\, u^{-1}$ and 
$\pi_2(a_2)= u^2\, (1\ot a_2)\, u^{-2}$ have free images in $(\Tt_B ( E),\rho\circ\mathfrak{E})$. 
Thus they generate a $C(X)$-linear monomorphism 
$\pi: A\hookrightarrow \Tt_B ( E)$ extending each $\pi_k$ and 
satisfying $\rho\circ\mathfrak{E}\circ\pi=\phi$ 
(Lemma~3.1 and Proposition~3.2 of \cite{DySh}, or \cite{BlDy}). 
Above Lemma~\ref{lemToeplitz} entails that $A$ is a continuous $C(X)$-algebra with fibre at $x\in X$ its image in $\Tt_{B_x}(E_x)$, i.e. the reduced free product $((A_1)_x,(\phi_1)_x)\ast\,((A_2)_x,(\phi_2)_x)$. 
\qed

\begin{rem} 
The existence of an embedding of $A_1$ in $C(X; \Od)$ 
cannot give a direct proof of Theorem \ref{redfree} 
since there is no $C(X)$-linear Hahn-Banach theorem (\cite[4.2]{bla3}). 
\end{rem}

\begin{cor} 
Any separable unital continuous $C(X)$-algebra $A$ admits 
a $C(X)$-linear unital embedding into a unital continuous field 
$\widetilde{A}$ with simple fibres. 
\end{cor} 
\pf Let $\phi: A\to C(X)$ be a $C(X)$-linear unital map such that 
each induced state $\phi_x: A_x\to\C$ is faithful. 
Then, the reduced free product
$$(\widetilde{A},\Phi)=(A\otimes\C^2;\phi\otimes tr_2)\ast_{C(X)}
(C(X)\otimes\C^3;\mathrm{id}\otimes tr_3)$$ is continuous 
by Proposition~\ref{redfree}, and it has simple fibres (\cite{Av}, \cite{Ba}).
\qed

\begin{cor}\label{cor 4.8} 
Let $X$ be a second countable \textsl{perfect} compact space and 
$A_1$ a unital separable continuous $C(X)$-algebra. 
Then the following assertions are equivalent.
\begin{itemize} 
\item[$\alpha )$] The \cst-algebra $A_1$ is exact. 
\item[$\beta )$] For all unital separable continuous $C(X)$-algebra $A_2$ and 
all continuous fields of faithful states $\phi_1$ and $\phi_2$ on $A_1$ and $A_2$, 
the reduced amalgamated free product 
$(A,\phi)=(A_1,\phi_1)\mathop{\ast}\limits_{C(X)} (A_2,\phi_2)$ is 
a continuous $C(X)$-algebra 
with fibres $(A_x,\phi_x)=((A_1)_x,(\phi_1)_x) \ast \,((A_2)_x,(\phi_2)_x)$. 
\end{itemize}
\end{cor}
\pf 
We only need to prove the implication $\beta)\Rightarrow\alpha)$ since the reverse implication has already been proved in Theorem~\ref{redfree}. 

Now, if a pair $(A_2, \phi_2)$ satisfies the hypotheses of $\beta )$ and 
we define the $C(X)$-algebra $B:=A_1\otimes_{C(X)} A_2$, 
the $C(X)$-linear projection $\rho=\phi_1\otimes\phi_2: B\to C(X)$ and 
the Hilbert $B$-module $E=L^2(B, \rho)\otimes_{C(X)} B$, 
then we have a $C(X)$-linear isomorphism 
$A\rtimes_\alpha\mathbb{N}\cong \mathcal{T}_B(E\oplus B)$ (Step 2 of Lemma \ref{lemToeplitz}). 
Hence, the Toeplitz $C(X)$-algebra $\mathcal{T}_B(E\oplus B)$ is continuous 
since the group $\mathbb{Z}$ is amenable (see e.g. [30]). 
And so, the amalgamated tensor product $A_1\otimes_{C(X)} A_2$ is a continuous $C(X)$-algebra for any unital separable continuous $C(X)$-algebra $A_2$ (Lemma \ref{lemToeplitz}). 
But this implies the exactness of the \cst-algebra $A_1$ 
if the metrizable space $X$ is perfect (\cite[Theorem 1.1]{BdWa}). 
\qed

\begin{rem} 
Corollary \ref{cor 4.8} does not always hold if the space $X$ is not perfect. 
For instance, if $X$ is reduced to a point, then the reduced amalgamated free product of $A_1$ and $A_2$ is always continuous. 
\end{rem}

\noindent
\email{Etienne.Blanchard@math.jussieu.fr}\\
\address{IMJ,
175, rue du Chevaleret, F--75013 Paris}

\end{document}